\documentclass[a4paper,11pt,twoside]{article}
\usepackage{amsfonts}
\usepackage{latexsym,amssymb}
\usepackage{amsmath}

\newtheorem{example}{Example}[section]
\newtheorem{defn}[example]{Definition}
\newtheorem{prop}[example]{Proposition}
\newtheorem{thm}[example]{Theorem}
\newtheorem{cor}[example]{Corollary}

\newenvironment{pf}{{\bf Proof:} }{$\Box$
\mbox{}}

\input xypic
\input xy
\xyoption{all}
\input{tcilatex}
\begin{document}

\author{\.{I}brahim \.{I}lker Ak\c{c}a and Ummahan Ege Arslan }
\title{Categorification of Algebras: 2-Algebras}
\date{}
\maketitle

\begin{abstract}
This paper introduces a categorification of $k$-algebras called $2$%
-algebras, where $k$ is a commutative ring. We define the $2$-algebras as a $%
2$-category with single object in which collections of all $1$-morphisms and
all $2$-morphisms are $k$-algebras. It is shown that the category of $2$%
-algebras is equivalent to the category of crossed modules in commutative $k$%
-algebras. Also we define the notion of homotopy for 2-algebras and we
explore the relations of crossed module homotopy and 2-algebra homotopy.
\end{abstract}

\footnotetext[1]{%
}

\section*{Introduction}

The term \textquotedblleft categorification\textquotedblright \ coined by
Louis Crane refers to the process of replacing set theoretic concepts by
category-theoretic analogues in mathematics. A categorified
version of a group is a 2-group. Internal categories in the category of groups are
exactly the same as 2-groups. The Brown-Spencer theorem \cite{BROWN-SPENCER}
thus constructs the associated 2-group of a crossed module given by
Whitehead \cite{WHITEHEAD} to define an algebraic model for a
\textquotedblleft(connected) homotopy 2-type\textquotedblright. The fact
that the composition in the internal category must be a group homomorphism
implies that the \textquotedblleft interchange law\textquotedblright \ must
hold. This equation is in fact equivalent via the Brown-Spencer result to
the Peiffer identity.

We will be concerned in this paper exclusively with categorification of
algebras. We will obtain analogous results in (commutative) algebras with
regard to Porter's work \cite{TIM}. He states that there is an equivalence
of categories between the category of internal categories in the category of $k$-algebras and
the category of crossed modules of commutative $k$-algebras. Since the internal category in the
category of $k$-algebras \ is a categorification of $k$-algebras, this
internal category will be called as \textquotedblleft strict
2-algebra\textquotedblright\ in this work. We define the strict 2-algebra by means of $2$-module being a category in the category of modules as
a 2-category with single object in which collections of 1-morphisms and
2-morphisms are $k$-algebras and we denote the category of strict 2-algebras
by \textbf{2Alg }. Given a group $G$, it is known that automorphisms of $G$
yield a $2$-group. Analogous result in algebras can be given that
multiplications of $C$ yield a strict $2$-algebra where $C$ is an $R$%
-algebra and $R$ is a $k$-algebra.

A crossed module %\cite{PORTER}
$\mathcal{A}=(\partial :C\longrightarrow R)$ of commutative algebras is
given by an algebra morphism $\partial :C\longrightarrow R$ together with an
action $\cdot $ of $R$ on $C$ such that the relations below hold for each $%
r\in R$ and each $c,c^{^{\prime }}\in C,$%
\begin{equation*}
\begin{array}{ccc}
\partial (r\cdot c) & = & r\partial (c) \\
\partial (c)\cdot c^{^{\prime }} & = & cc^{^{\prime }}.%
\end{array}%
\end{equation*}%
Group crossed modules were firstly introduced by Whitehead in \cite{WH},\cite%
{WHT}. They are algebraic models for homotopy 2-types, in the sense that
\cite{BAUES},\cite{LODAY} the homotopy category of the model category \cite%
{BROWN-GOLAS},\cite{CABELLO} of group crossed modules is equivalent to the
homotopy category of the model category \cite{ELVIRA} of pointed 2-types:
pointed connected spaces whose homotopy groups $\pi _{i}$ vanish, if $i\geq
3 $. The homotopy relation between crossed module maps $\mathcal{A}%
\longrightarrow \mathcal{A}^{^{\prime }}$ was given by Whitehead in \cite%
{WHT}, in the contex of \textquotedblleft homotopy
systems\textquotedblright\ called free crossed complexes.

In \cite{kim} it is addressed the homotopy theory of maps between crossed modules of commutative algebras. It is proven that if $\mathcal{A}$
and $\mathcal{A}^{^{\prime }}$ are crossed modules of algebras without
any restriction on $\mathcal{A}$ and $\mathcal{A}^{^{\prime }}$ then the crossed module maps $\mathcal{A}\longrightarrow \mathcal{A}^{^{\prime }}$
and their homotopies give a groupoid.

In this paper we show that the category of strict 2-algebras is equivalent
to the category of crossed modules in commutative algebras. In \cite{E.K}, it is given an equivalence between the category of crossed modules in associative algebras and the category of strict associative $2$-algebras defined by means of $2$-vector space. Also we define
the notion of homotopy for 2-algebras. This definition is essentially a
special case of 2-natural transformation due to Gray in \cite{Gray}. And we
explore the relations between the crossed module homotopies and 2-algebra
homotopies. Similar results are given \cite{icen} by \.{I}\c{c}en for
2-groupoids.

%\begin{equation*}
%\xymatrix{{\textbf{XMod}_{k}}\ar@{<->}[dddr]_{\cite{NIZAR}}\ar@<0.5ex>[dr]^{%
%\Gamma}&&\textbf{SimpAlg$_{\leq
%1}$}\ar[dddl]^{\cite{LODAY}}\ar@{<->}[ll]_{\cite{ARVASI}}\ar@<0.5ex>[dl]^{%
%\Omega}\\ &\textbf{2Alg}\ar@<0.5ex>[ur]^{\chi}
%\ar@<0.5ex>[dd]^{\Lambda}\ar@<0.5ex>[ul]^{\Psi}&\\ \\
%&\textbf{Cat}^{1}\ar@<0.5ex>[uu]^{\Delta}\ar[uuur]&&}
%\end{equation*}

\section{Internal Categories and 2-categories}

%{\bf Internal Categories} \vskip12 pt

We begin by recalling internal categories as well as $2$-categories.
Ehresmann defined internal categories in \cite{EHR}, and
by now they are an important part of category theory \cite{BOR}.

\subsection{Internal categories}

\begin{defn}
\label{1.1} Let $\mathbf{C}$ be any category. An internal category in $\mathbf{C}$, say $%
\mathbf{A}$, consists of:

$\ {\small \blacklozenge } $ an object of objects $A_{0}\in \mathbf{C}$

$\ {\small \blacklozenge } $ an object of morphisms $A_{1}\in \mathbf{C}$,

\noindent together with

$\ {\small \blacklozenge } $ source and target morphisms $s,t:A_{1}\longrightarrow A_{0}%
\mathbf{,}$

$\ {\small \blacklozenge } $ an identity-assigning morphism $e:A_{0}\longrightarrow A_{1},$

$\ {\small \blacklozenge } $ a composition morphism $\circ :A_{1}\times
_{A_{0}}A_{1}\longrightarrow A_{1}$ such that the following diagrams
commute, expressing the usual category laws:

$\ {\small \blacklozenge } $ laws specifying the source and target of identity morphisms:

\begin{equation*}
\xymatrix { A_0 \ar[drr]_{1_{A_0}}\ar[rr]^{e} & &A_1\ar[d]^{s} & \\ \ & &
A_0& } \xymatrix { A_0 \ar[drr]_{1_{A_0}}\ar[rr]^{e} & &A_1\ar[d]^{t} & \\ \
& & A_0& }
\end{equation*}
$\ {\small \blacklozenge } $ laws specifying the source and target of composite morphisms:

\begin{equation*}
\xymatrix{ &&A_1 \times _{A_{0}}A_1\ar[dd]_{\rho_{1}} \ar[rr]^{\circ } &&
A_1 \ar[dd]^{s} && \\ \\ &&A_1 \ar[rr]_{s} && A_0 && }
\end{equation*}%
\begin{equation*}
\xymatrix{ &&A_1 \times_{A_{0}}A_1\ar[dd]_{\rho_{2}} \ar[rr]^{\circ } && A_1
\ar[dd]^{t} && \\ \\ &&A_1 \ar[rr]_{t} && A_0 && }
\end{equation*}%
$\ {\small \blacklozenge } $ the associative law for composition of morphisms:
\begin{equation*}
\xymatrix{ &&A_1 \times _{A_{0}}A_1 \times _{A_{0}}A_1 \ar[dd]_{\rho_{2}}
\ar[rr]^{\circ } && A_1 \times _{A_{0}}A_1 \ar[dd]^{t} && \\ \\ &&A_1 \times
_{A_{0}}A_1 \ar[rr]_{t} && A_0 && }
\end{equation*}%
$\ {\small \blacklozenge } $ the left and right unit laws for composition of morphisms:
\begin{equation*}
\xymatrix { A_0 \times_{A_0} A_1 \ar[r]^{ e \times_{A_0} 1_{A_1} }
\ar[dr]_{\rho_2} & A_1 \times_{A_0} A_1 \ar[d]^{\circ } & A_1 \times_{A_0}
A_0 \ar[l]_{1_{A_1} \times_{A_0} e} \ar[dl]^{\rho_1} & \\ \ & A_1 & & }
\end{equation*}%
The pullback $A_{1}\times _{A_{0}}A_{1}$ is defined via the square:
\begin{equation*}
\xymatrix { A_1 \times _{A_{0}}A_1 \ar[r]^{\rho_{2}} \ar[d]_{\rho_{1}} & A_1
\ar[d]^{s} & \\ \ A_1 \ar[r]_{t} & A_0. & }
\end{equation*}%
We denote this internal category with $A=(A_{0},A_{1},s,t,e,\circ )$.
\end{defn}

\begin{defn}
Let $\mathbf{C}$ be a category. Given internal categories $A$ and $A^{\prime
}$ in $\mathbf{C}$, an \textbf{internal functor} \ between them, say $%
F:A\longrightarrow A^{\prime }$, consists of

$\ {\small \blacklozenge } $ a morphism $F_{0}:A_{0}\longrightarrow A_{0}^{\prime },$

$\ {\small \blacklozenge } $ a morphism $F_{1}:A_{1}\longrightarrow A_{1}^{\prime }$

\noindent such that the following diagrams commute, corresponding to the usual laws
satisfied by a functor:

$\ {\small \blacklozenge } $ preservation of source and target:%
\begin{equation*}
\xymatrix { A_1 \ar[d]_{F_1}\ar[r]^{s} & A_0\ar[d]^{F_0} \\ \ A_{1}^{\prime
}\ar[r]_{s^{\prime }}& A_{0}^{\prime } } \xymatrix { A_1
\ar[d]_{F_1}\ar[r]^{t} & A_0\ar[d]^{F_0} \\ \ A_{1}^{\prime
}\ar[r]_{t^{\prime }}& A_{0}^{\prime } }
\end{equation*}

$\ {\small \blacklozenge } $ preservation of identity morphisms:%
\begin{equation*}
\xymatrix { A_0 \ar[d]_{F_0}\ar[r]^{e} & A_1\ar[d]^{F_1} \\ \ A_{0}^{\prime
}\ar[r]_{e^{\prime }}& A_{1}^{\prime } }
\end{equation*}

$\ {\small \blacklozenge } $ preservation of composite morphisms:%
\begin{equation*}
\xymatrix{ &&A_1 \times _{A_{0}}A_1 \ar[dd]_{\circ} \ar[rr]^{F_{1} \times
_{A_{0}}F_{1}} && A_{1}^{\prime } \times _{A_{0}^{\prime }}A_{1}^{\prime }
\ar[dd]^{{\circ}^{\prime}} && \\ \\ &&A_1 \ar[rr]_{F_1} && A_{1}^{\prime} &&
}
\end{equation*}
\end{defn}

Given two internal functors $F:A\longrightarrow A^{\prime }$ and $%
G:A^{\prime }\longrightarrow A^{\prime \prime }$ in some category $\mathbf{C}
$, we define their composite $FG:A\longrightarrow A^{\prime \prime }$ by
taking $(FG)_{0}=F_{0}G_{0}$ and $(FG)_{1}=F_{1}G_{1}$. Similarly, we define
the identity internal functor in $\mathbf{C}$, $1_{A}:A\longrightarrow A$ by
taking $(1_{A})_{0}=1_{A_{0}}$ and $(1_{A})_{1}=1_{A_{1}}$.

\begin{defn}
Let $\mathbf{C}$ be a category. Given two internal functors $%
F,G:A\longrightarrow A^{\prime }$ in $\mathbf{C}$, an \textbf{internal
natural transformation} in $\mathbf{C}$ between them, say $\theta
:F\Longrightarrow G$, is a morphism $\theta :A_{0}\longrightarrow
A_{1}^{\prime }$ for which the following diagrams commute, expressing the
usual laws satisfied by a natural transformation:

$\ {\small \blacklozenge } $ laws specifying the source and target of a natural transformation:%
\begin{equation*}
\xymatrix { A_0 \ar[drr]_{F_0}\ar[rr]^{\theta} &
&A_{1}^{\prime}\ar[d]^{{s}^{\prime}} & \\ \ & & A_{0}^{\prime}& } \xymatrix
{ A_{0}^{\prime} \ar[drr]_{G_0}\ar[rr]^{\theta} &
&A_{1}^{\prime}\ar[d]^{{t}^{\prime}} & \\ \ & & A_0& }
\end{equation*}

$\ {\small \blacklozenge } $ the commutative square law:%
\begin{equation*}
\xymatrix{ &&A_1 \ar[dd]_{\vartriangle(F\times t\theta)}
\ar[rr]^{\vartriangle(s\theta \times G)} && A_{1}^{\prime } \times
_{A_{0}^{\prime }}A_{1}^{\prime } \ar[dd]^{{\circ}^{\prime}} && \\ \\
&&A_{1}^{\prime } \times _{A_{0}^{\prime }}A_{1}^{\prime
}\ar[rr]_{{\circ}^{\prime}} && A_{1}^{\prime} && }
\end{equation*}
\end{defn}

Given an internal functor $F:A\longrightarrow A^{\prime }$ in $\mathbf{C}$,
the identity internal natural transformation $1_{F}:F\Longrightarrow F$ in $%
C $ is given by $1_{F}=F_{0}e$.

\subsection{2-categories}

\begin{definition}
A $2$-category $\mathcal{G}$ consists of a class of objects $G_{0}$ and for
any pair of objects $(A,B)$ a small category of morphisms $\mathcal{G}(A,B)$%
-with objects $G_{1}(A,B)$ and morphisms $G_{2}(A,B)$-, along with
composition functors
\begin{equation*}
~%
\begin{array}{ccccc}
\bullet  & : & \mathcal{G}(A,B)\times \mathcal{G}(B,C) & \longrightarrow  &
\mathcal{G}(A,C)%
\end{array}%
\end{equation*}%
for every triple $(A,B,C)$ of objects and identity functors from the
terminal category to $\mathcal{G}(A,A)$
\begin{equation*}
iA:1\longrightarrow \mathcal{G}(A,A)
\end{equation*}%
for all objects $A$ such that $\bullet $ \ is\ associative\ and
\begin{equation*}
F\bullet i_B=F=i_A\bullet F\quad \text{as well as}\quad \vartheta \bullet I_{i_B}
=\vartheta =I_{i_A}\bullet \vartheta
\end{equation*}%
hold for all $F\in G_{1}(A,B)$ and $\vartheta \in G_{2}(A,B)$ where source and target morphisms are defined by %
\begin{equation*}
A\overset{F}{\longrightarrow }B
\end{equation*}%
\begin{equation*}
\begin{array}{ccc}
s:G_{1}(A,B) & \longrightarrow & G_{0} \\
F & \longmapsto & s(F)=A%
\end{array}%
\end{equation*}%
\begin{equation*}
\begin{array}{ccc}
t:G_{1}(A,B) & \longrightarrow & G_{0} \\
F & \longmapsto & t(F)=B%
\end{array}%
\end{equation*}%
for $F\in G_{1}(A,B)$ and %
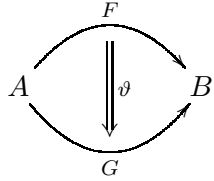
\begin{equation*}
\xymatrix{ A \ar@/^2pc/[rr]_{\quad}^{F}="1" \ar@/_2pc/[rr]_{G}="2" && B
\ar@{}"1";"2"|(.2){\,}="7" \ar@{}"1";"2"|(.8){\,}="8" \ar@{=>}"7"
;"8"^{\vartheta} }
\end{equation*}%

\begin{equation*}
\begin{array}{ccc}
s:G_{2}(A,B) & \longrightarrow & G_{1} \\
\vartheta & \longmapsto & s(\vartheta )=F%
\end{array}%
\end{equation*}%
\begin{equation*}
\begin{array}{ccc}
t:G_{2}(A,B) & \longrightarrow & G_{0} \\
\vartheta & \longmapsto & t(\vartheta )=G%
\end{array}%
\end{equation*}%
for $\vartheta :F\longrightarrow G\in G_{2}(A,B)$.
For all pairs of objects $(A,B)$ elements of $G_{1}(A,B)$ are called $1$%
-morphisms or $1$-cells of $\mathcal{G}$ and elements of $G_{2}(A,B)$ are
called $2$-morphisms or $2$-cells of $\mathcal{G}$. We write $G_{1}$ and $%
G_{2}$ for the classes of all $1$-morphisms and $2$-morphisms respectively.

There are two ways of composing $2$-morphisms: using the composition $%
\circ $ inside the categories $\mathcal{G}(A,B)$, called vertical
composition, and using the morphism level of the functor $\bullet $, called
horizontal composition. These compositions must be satisfy the following
equation: for $\alpha ,\alpha ^{\prime }\in G_{2}(A,B)$ with $t(\alpha
)=s(\alpha ^{\prime })$ and $\gamma ,\gamma ^{\prime }\in G_{2}(B,C)$ with $%
t(\gamma )=s(\gamma ^{\prime })$

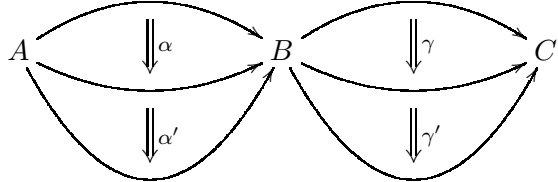
\begin{equation*}
\xymatrix{ A \ar@/^1.6pc/[rrr]_{\quad}^{}="1" \ar@/_1.2pc/[rrr]_{ }="2"
\ar@/_4pc/[rrr]_{}="9" &&& B \ar@{}"1";"2"|(.2){\,}="7"
\ar@{}"1";"2"|(.8){\,}="8" \ar@{=>}"7" ;"8"^{\alpha}
\ar@{}"2";"9"|(.2){\,}="10" \ar@{}"2";"9"|(.8){\,}="11" \ar@{=>}"10"
;"11"^{\alpha^{\prime }} \ar@/^1.6pc/[rrr]_{\quad}^{}="12"
\ar@/_1.2pc/[rrr]_{}="13"\ar@/_4pc/[rrr]_{}="14"
\ar@{}"12";"13"|(.2){\,}="15" \ar@{}"12";"13"|(.8){\,}="16" \ar@{=>}"15"
;"16"^{\gamma} \ar@{}"13";"14"|(.2){\,}="17" \ar@{}"13";"14"|(.8){\,}="18"
\ar@{=>}"17" ;"18"^{\gamma ^{\prime }} &&& C }
\end{equation*}%
\begin{equation*}
\left( \alpha \circ \alpha ^{\prime }\right) \bullet \left( \gamma \circ
\gamma ^{\prime }\right) =\left( \alpha \bullet \gamma \right) \circ \left(
\alpha ^{\prime }\bullet \gamma ^{\prime }\right)
\end{equation*}%
which is called \textquotedblleft interchange law\textquotedblright .
\end{definition}

\section{Constructions of Two-Algebras}

In this section we will construct $2$-algebras by categorification. We can categorify the notion of an algebra by replacing the equational laws by
isomorphisms satisfying extra structure and properties we expect. In \cite%
{CRANS} Baez and Crans introduce the Lie 2-algebra by means of the
concept of $2$-vector space defined as an internal category in the category
of vector spaces by them. Obviously we get a new notion of \textquotedblleft
2-module\textquotedblright which can be considered as an internal category
in the category of modules and we categorify the notion of an algebra.

\subsection{2-Modules}

A categorified module or \textquotedblleft 2-module\textquotedblright
should be a category with structure analogous to that of a $k-$module, with
functors replacing the usual $k-$module operations. Here we instead define a
2-module to be an internal category in a category of $k-$modules \textbf{Mod
}. Since the main component part of a $k-$algebra is a $k-$module, a 2-algebra
will have an underlying 2-module of this sort. In this section we thus first
define a category of these 2-modules.

In the rest of this paper, the terms a module and an algebra will always
refer to a $k-$module and a $k-$algebra.

\begin{defn}
A 2-module is an internal category in \textbf{Mod }.
\end{defn}

Thus, a 2-module $M$ is a category with a module of objects $M${}$_{0}$ and
a module of morphisms $M_{1}$, such that the source and target maps $%
s,t:M_{1}\longrightarrow M_{0}$, the identity assigning map $%
e:M_{0}\longrightarrow M_{1}$, and the composition map $\circ :M_{1}\times
_{M_{0}}M_{1}\longrightarrow M_{1}$ are all module morphisms. We
write a morphism as $a:x\longrightarrow y$ when $s(a)=x$ and $t(a)=y$, and
sometimes we write $e(x)$ as $1_{x}$.

The following proposition is given for the \textbf{Vect } vector space category in
\cite{CRANS}. But we rewrite this proposition for \textbf{Mod }.

\begin{prop}
\label{2.2} It is defined a $2$-module by specifying the modules $M_{0}$ and $M_{1}$ along with the source, target and identity
module morphisms and the composition morphism $\circ $, satisfying the conditions of Definition \ref{1.1}. The composition map is uniquely determined by
\begin{equation*}
\begin{array}{lcll}
\circ : & M_{1}\times _{M_{0}}M_{1} & \longrightarrow & M_{1} \\
& (a,b) & \longmapsto & \circ (a,b)=a\circ b=a+b-(es)(b).%
\end{array}%
\end{equation*}
\end{prop}

\begin{pf}
First given modules $M_{0}$, $M_{1}$ and module morphisms $%
s,t:M_{1}\longrightarrow M_{0}$ and $e:M_{0}\longrightarrow M_{1}$, we will
define a composition operation that satisfies the laws in the definition of
internal category, obtaining a 2-module.

Given $a$, $b\in M_{1}$ such that $t(a)=s(b)\,$, i.e.%
\begin{equation*}
a:x\longrightarrow y\text{ and }b:y\longrightarrow z
\end{equation*}%
we define their composite $\circ $ by%
\begin{equation*}
\begin{array}{lcll}
\circ : & M_{1}\times _{M_{0}}M_{1} & \longrightarrow & M_{1} \\
& (a,b) & \longmapsto & \circ (a,b)=a\circ b=a+b-(es)(b).%
\end{array}%
\end{equation*}%
We will show that with this composition $\circ $ the diagrams of
the definition of internal category commute. The triangles specifying the source
and target of the identity-assigning morphism do not involve composition.
The second pair of diagrams commute since
\begin{equation*}
\begin{array}{lll}
s(a\circ b) & = & s(a+b-(es)(b)) \\
& = & s(a)+s(b)-(se)(s(b)) \\
& = & s(a)+s(b)-s(b) \\
& = & s(a)=x%
\end{array}%
\end{equation*}%
and since $t(a)=s(b),$%
\begin{equation*}
\begin{array}{lll}
t(a\circ b) & = & t(a+b-(es)(b)) \\
& = & t(a)+t(b)-(te)(s(b)) \\
& = & t(a)+t(b)-s(b) \\
& = & t(b)=z.%
\end{array}%
\end{equation*}%
Since module operation is associative, the associative law holds for composition. The left and right unit laws are satisfied since given $%
a:x\longrightarrow y,$%
\begin{equation*}
\begin{array}{lll}
e(x)\circ a & = & e(x)+a-(es)(a) \\
& = & e(x)+a-e(x) \\
& = & a%
\end{array}%
\end{equation*}%
and%
\begin{equation*}
\begin{array}{lll}
a\circ e(y) & = & a+e(y)-(es)(e(y)) \\
& = & a+e(y)-e(y) \\
& = & a.%
\end{array}%
\end{equation*}%
We thus have a 2-module.

Given a 2-module $M$, we show that its composition must be
defined by the formula given above. Let $(a,g)$ and $(a^{\prime
},g^{\prime })$ be composable pairs of morphisms in $M_{1}$, i.e.%
\begin{equation*}
a:x\longrightarrow y\text{ and }b:y\longrightarrow z
\end{equation*}%
and%
\begin{equation*}
a^{\prime }:x^{\prime }\longrightarrow y^{\prime }\text{ and }b^{\prime
}:y^{\prime }\longrightarrow z^{\prime }.
\end{equation*}%
Since the source and target maps are module morphisms, $(a+a^{\prime
},b+b^{\prime })$ also forms a composable pair, and since that the
composition is module morphism%
\begin{equation*}
(a+a^{\prime })\circ (b+b^{\prime })=a\circ b+a^{\prime }\circ b^{\prime }.
\end{equation*}

Then if $(a,b)$ is a composable pair, i.e, $t(a)=s(b),$ we have%
\begin{equation*}
\begin{array}{lll}
a\circ b & = & (a+1_{M_{1}})\circ (1_{M_{1}}+b) \\
& = & (a+e(s(b)-s(b)))\circ (e(s(b)-s(b))+b) \\
& = & (a-e(s(b))+e(s(b)))\circ (e(s(b))-e(s(b))+b) \\
& = & (a\circ e(s(b)))+(-e(s(b))+e(s(b)))\circ (-e(s(b))+b) \\
& = & a\circ e(s(b))+(-e(s(b))\circ (-e(s(b))))+(e(s(b))\circ b) \\
& = & a-e(s(b))+b \\
& = & a+b-e(s(b)).%
\end{array}%
\end{equation*}%
This show that we can define $\circ $ by%
\begin{equation*}
\begin{array}{lcll}
\circ : & M_{1}\times _{M_{0}}M_{1} & \longrightarrow & M_{1} \\
& (a,b) & \longmapsto & \circ (a,b)=a\circ b=a+b-e(s(b)).%
\end{array}%
\end{equation*}
\end{pf}

\begin{cor}
For $b\in \ker s,$ we have
\begin{equation*}
\begin{array}{lll}
a\circ b & = & a+b-(es)(b) \\
& = & a+b.%
\end{array}%
\end{equation*}
\end{cor}

\begin{defn}
Let $M$ and $N$ be 2-modules, a 2-module functor $F:M\longrightarrow N$ is
an internal functor in $\mathbf{Mod}$ from $M$ to $N$. 2-modules and
2-module functors between them is called the category of 2-modules denoted by $%
\mathbf{2Mod}$.
\end{defn}

After we get the definition of a 2-module, we define the definition of a categorified algebra which is main concept of this paper.

\subsection{Two-algebras}

\begin{defn}
A weak 2-algebra consists of

$\ {\small \blacklozenge }$ a 2-module $A$ equipped with a functor $\bullet :A\times
A\longrightarrow A$, which is defined by $(x,y)\mapsto x\bullet y$ and
bilinear on objects and defined by $(f,g)\mapsto f\bullet g$ on morphisms
satisfying interchange law, i.e.,
\begin{equation*}
(f_{1}\bullet g_{1})\circ (f_{2}\bullet g_{2})=(f_{1}\circ f_{2})\bullet
(g_{1}\circ g_{2})
\end{equation*}

$\ {\small \blacklozenge }$ $k-$bilinear natural isomorphisms%
\begin{equation*}
\alpha _{x,y,z}:(x\bullet y)\bullet z\longrightarrow x\bullet (y\bullet z)
\end{equation*}%
\begin{equation*}
l_{x}:1\bullet x\longrightarrow x
\end{equation*}%
\begin{equation*}
r_{x}:x\bullet 1\longrightarrow x
\end{equation*}

such that the following diagrams commute for all objects $w,x,y,z\in A_{0}.$

\begin{equation*}
\xymatrix{ ((w\bullet x)\bullet y)\bullet z
\ar[d]_{{{\alpha}_{w,x,y}}\bullet {1_z}} \ar[r]^{{\alpha}_{w\bullet x,y,z}}
& (w\bullet x)\bullet(y\bullet z) \ar[dr]^{{\alpha}_{w,x,y\bullet z}} \\
(w\bullet (x\bullet y))\bullet z \ar[r]_{{\alpha}_{w,x\bullet y,z}} &
w\bullet ((x\bullet y)\bullet z) \ar[r]_{{1_w}\bullet {{\alpha}_{x,y,z}}} &
w\bullet (x\bullet (y\bullet z)) }
\end{equation*}

\begin{equation*}
\xymatrix{ (x\bullet 1)\bullet y \ar[dr]_{{r_x}\bullet {1_y}}
\ar[r]^{{\alpha}_{x,1,y}} & x\bullet (1\bullet y) \ar[d]^{{1_x}\bullet
{l_y}} \\ & x\bullet y }
\end{equation*}

A strict 2-algebra is the special case where $\alpha _{x,y,z}$, $l_{x}$, $%
r_{x}$ are all identity morphisms. In this case we have
\begin{equation*}
(x\bullet y)\bullet z=x\bullet (y\bullet z)
\end{equation*}
\begin{equation*}
1\bullet x=x, x\bullet 1=x
\end{equation*}%
Strict 2-algebra is called commutative strict 2-algebra if $x\bullet
y=y\bullet x$ for all objects $x,y\in A_{0}$ and $f\bullet g=g\bullet f$ for
all morphisms $f,g\in A_{1}$.
\end{defn}

In the rest of this paper, the term 2-algebra will always refer to a
commutative strict 2-algebra. A homomorphism between 2-algebras should
preserve both the 2-module structure and the $\bullet $ functor.

\begin{defn}
Given 2-algebras $A$ and $A^{\prime }$, a homomorphism%
\begin{equation*}
F:A\longrightarrow A^{\prime }
\end{equation*}%
consists of

$\ {\small \blacklozenge }$ a linear functor $F$ from the underlying 2-module of $A$ to that of
$A^{\prime }$, and

$\ {\small \blacklozenge }$ a bilinear natural transformation%
\begin{equation*}
F_{2}(x,y):F_{0}(x)\bullet F_{0}(y)\longrightarrow F_{0}(x\bullet y)
\end{equation*}

$\ {\small \blacklozenge }$ an isomorphism $F:1^{\prime }\longrightarrow F_{0}(1)$ where $1$ is
the identity object of $A$ and $1^{\prime }$ is the identity object of $%
A^{\prime }$,

such that the following diagrams commute for $x,y,z\in A_{0}$,

\begin{equation*}
\xymatrix{ (F(x)\bullet F(y))\bullet F(z) \ar[d]_{{\alpha}_{F(x),F(y),F(z)}}
\ar[r]^{{F_2}\bullet 1} & F(x\bullet y)\bullet F(z) \ar[r]^{F_2} &
F((x\bullet y)\bullet z) \ar[d]^{F({\alpha}_{x,y,z})} \\ F(x)\bullet
(F(y)\bullet F(z)) \ar[r]_{1\bullet {F_2}} & F(x)\bullet F(y\bullet z)
\ar[r]_{F_2} & F(x\bullet (y\bullet z)). }
\end{equation*}

\begin{equation*}
\xymatrix{ 1^{\prime }\bullet F(x) \ar[d]_{{F_0}\bullet 1}
\ar[r]^{{l^{\prime }}_{F(x)}} & F(x) \\ F(1)\bullet F(x) \ar[r]_{F_2} &
F(1\bullet x) \ar[u]_{F(l_x)}.}
\end{equation*}

\begin{equation*}
\xymatrix{ F(x)\bullet 1^{\prime } \ar[d]_{1\bullet {F_0}}
\ar[r]^{{r^{\prime }}_{F(x)}} & F(x) \\ F(x)\bullet F(1) \ar[r]_{F_2} &
F(x\bullet 1) \ar[u]_{F(r_x)}.}
\end{equation*}
\end{defn}

\begin{defn}
2-algebras and homomorphisms between them give the category of 2-algebras
denoted by \textbf{2Alg }.
\end{defn}

Therefore if $A=(A_{0},A_{1},s,t,e,\circ ,\bullet )$ is a 2-algebra, $A_{0}$
and $A_{1}$ are algebras with this $\bullet $ bilinear functor. Thus we can
take that 2-algebra is a 2-category with a single object say $\ast $, and $%
A_{0}$ collections of its 1-morphisms and $A_{1}$ collections of its
2-morphisms are algebras with identity.

\subsection{Multiplication Algebras yield a 2-algebra}

In \cite{N} Norrie developed Lue's work, \cite{L} and introduced the notion
of an actor of crossed modules of groups where it is shown to be the
analogue of the automorphism group of a group. In the
category of commutative algebras the appropriate replacement for
automorphism groups is the multiplication algebra $\mathcal{M}(C)$ of an
algebra $C$ which is defined by MacLane \cite{MACLANE}.

Let $C$ be an associative (not necessarily unitary or commutative) $R$%
-algebra. We recall Mac Lane's construction of the $R$-algebra Bim$(C)$ of
bimultipliers of $C$ \cite{MACLANE}.

An element of Bim$(C)$ is a pair $(\gamma ,\delta )$ of $R$-linear mappings
from $C$ to $C$ such that
\begin{equation*}
\gamma (cc^{\prime })=\gamma (c)c^{\prime }
\end{equation*}%
\begin{equation*}
\delta (cc^{\prime })=c\delta \left( c^{\prime }\right)
\end{equation*}%
and
\begin{equation*}
c\gamma \left( c^{\prime }\right) =\delta (c)c^{\prime }.
\end{equation*}%
Bim$(C)$ has an obvious $R$-module structure and a product
\begin{equation*}
(\gamma ,\delta )(\gamma ^{\prime },\delta ^{\prime })=(\gamma \gamma
^{\prime },\delta ^{\prime }\delta ),
\end{equation*}%
the value of which is still in Bim$(C).$

Suppose that Ann$(C)=0$ or $C^{2}=C$. Then Bim$(C)$ acts on $C$ by
\begin{eqnarray*}
\text{Bim}(C)\times C &\rightarrow &C;\qquad ((\gamma ,\delta ),c)\mapsto
\gamma (c), \\
C\times \text{Bim}(C) &\rightarrow &C;\qquad (c,(\gamma ,\delta ))\mapsto
\delta (c)
\end{eqnarray*}%
and there is a%
\begin{equation*}
\begin{array}{llll}
\mu : & C & \longrightarrow & \text{Bim}(C) \\
& c & \longmapsto & (\gamma _{c},\delta _{c})%
\end{array}%
\end{equation*}%
with
\begin{equation*}
\gamma _{c}(x)=cx\qquad \text{and \qquad }\delta _{c}(x)=xc.
\end{equation*}

\textit{Commutative case}: we still assume Ann$(C)=0$ or $C^{2}=C.$ If $C$
is a commutative $R$-algebra and $(\gamma ,\delta )\in $ Bim$(C),$ then $%
\gamma =\delta .$ This is because for every $x\in C:$%
\begin{equation*}
\begin{array}{lll}
x\delta (c) & = & \delta (c)x=c\gamma (x)=\gamma (x)c \\
& = & \gamma (xc)=\gamma (cx)=\gamma (c)x=x\gamma (c).%
\end{array}%
\end{equation*}%
Thus Bim$(C)$ may be identified with the $R$-algebra $\mathcal{M}(C)$ of
multipliers of $C.$ Recall that a multiplier of $C$ is a linear mapping $%
\lambda :C\longrightarrow C$ such that for all $c,c^{\prime }\in C$%
\begin{equation*}
\lambda (cc^{\prime })=\lambda (c)c^{\prime }.
\end{equation*}%
Also $\mathcal{M}(C)$ is commutative as%
\begin{equation*}
\lambda ^{\prime }\lambda (xc)=\lambda ^{\prime }(\lambda (x)c)=\lambda
(x)\lambda ^{\prime }(c)=\lambda ^{\prime }(c)\lambda (x)=\lambda \lambda
^{\prime }(cx)=\lambda \lambda ^{\prime }(xc)
\end{equation*}%
for any $x\in C.$ Thus $\mathcal{M}(C)$ is the set of all multipliers $%
\lambda $ such that $\lambda \gamma =\gamma \lambda $ for every multiplier $%
\gamma .$

In \cite{PORTER} Porter states that automorphisms of a group $G$ yield a
2-group. The appropriate analogue of this result in algebra case can be
given. We claim that multiplications of an $R $-algebra $C$ give a 2-algebra
which is called a multiplication 2-algebra.

Let $k$ be a commutative ring, $R$ be a $k$-algebra with identity and $C$ be
a commutative $R$-algebra with $Ann(C)=0$ or $C^{2}=C$. Take $A_{0}=\mathcal{%
M}(C)$ and say 1-morphisms to the elements of $A_{0}$. We define the action
of $\mathcal{M}(C)$ on $C$ as follows:
\begin{equation*}
\begin{array}{ccc}
\mathcal{M}(C)\times C & \longrightarrow & C \\
(f,x) & \longmapsto & f\blacktriangleright x=f(x).%
\end{array}%
\end{equation*}%
Using the action of $\mathcal{M}(C)$ on $C$, we can form the semidirect
product
\begin{equation*}
C\rtimes \mathcal{M}(C)=\{(x,f)|x\in C,\text{ }f\in \mathcal{M}(C)\}
\end{equation*}%
with multiplication
\begin{equation*}
(x,f)(x^{\prime },f^{\prime })=(f\blacktriangleright x^{\prime }+f^{\prime
}\blacktriangleright x+x^{\prime }x,f^{\prime }f).
\end{equation*}%
Take $A_{1}=C\rtimes \mathcal{M}(C)$ and say 2-morphisms to the elements of $%
A_{1}$. Therefore we get the following diagram for $(x,f)\in C\rtimes
\mathcal{M}(C)$,

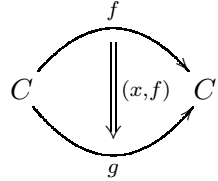
\begin{equation*}
\xymatrix{ C \ar@/^2pc/[rr]_{\quad}^{f}="1" \ar@/_2pc/[rr]_{g}="2" && C
\ar@{}"1";"2"|(.2){\,}="7" \ar@{}"1";"2"|(.8){\,}="8" \ar@{=>}"7"
;"8"^{(x,f)} }
\end{equation*}

and we define the source, target and identity assigning maps as follows;
\begin{equation*}
\begin{array}{ccccccccc}
s: & C\rtimes \mathcal{M}(C) & \longrightarrow & \mathcal{M}(C) &  & t: &
C\rtimes \mathcal{M}(C) & \longrightarrow & \mathcal{M}(C) \\
& (x,f) & \longmapsto & s(x,f)=f &  &  & (x,f) & \longmapsto & t(x,f)=M_{x}\cdot f%
\end{array}%
\end{equation*}%
and

\begin{equation*}
\begin{array}{cccc}
e: & \mathcal{M}(C) & \longrightarrow & C\rtimes \mathcal{M}(C) \\
& f & \longmapsto & e(f)=(0,f)%
\end{array}%
\end{equation*}%
where $M_{x}\cdot f$ is defined \ by $(M_{x}\cdot f)(u)=xu+f(u)$ for $u\in
C. $

There are two ways of composing 2-morphisms: vertical and horizontal
composition. Now we define this compositions.

For $(x,f),(y,f^{\prime })\in C\rtimes \mathcal{M}(C)$
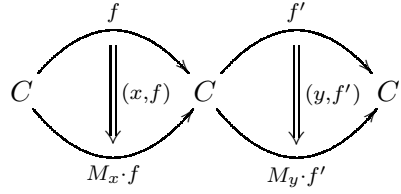
\begin{equation*}
\xymatrix{ C \ar@/^2pc/[rr]_{\quad}^{f}="1" \ar@/_2pc/[rr]_{M_{x}\cdot
f}="2" && C \ar@{}"1";"2"|(.2){\,}="7" \ar@{}"1";"2"|(.8){\,}="8"
\ar@{=>}"7" ;"8"^{(x,f)} \ar@/^2pc/[rr]_{\quad}^{f'}="9"
\ar@/_2pc/[rr]_{M_{y}\cdot f'}="10" \ar@{}"9";"10"|(.2){\,}="11"
\ar@{}"9";"10"|(.8){\,}="12" \ar@{=>}"11" ;"12"^{(y,f')} && C }
\end{equation*}%
the horizontal composition is defined by
\begin{equation*}
(x,f)\bullet (y,f^{\prime })=(f^{\prime }(x)+f(y)+xy,f^{\prime }f),
\end{equation*}%
thus we have
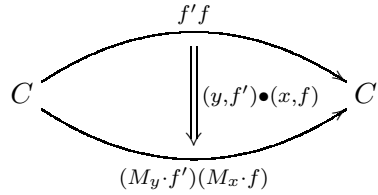
\begin{equation*}
\xymatrix{ C \ar@/^2pc/[rrrr]_{\quad}^{f'f}="1"
\ar@/_2pc/[rrrr]_{(M_{y}\cdot f')(M_{x}\cdot f)}="2" &&&& C
\ar@{}"1";"2"|(.2){\,}="7" \ar@{}"1";"2"|(.8){\,}="8" \ar@{=>}"7"
;"8"^{(y,f')\bullet (x,f)} }
\end{equation*}%
and
\begin{equation*}
\begin{array}{rcl}
t(f^{\prime }(x)+f(y)+xy,f^{\prime }f) & = & M_{f^{\prime }(x)+f(y)+xy}\cdot
f^{\prime }f \\
& = & (M_{y}\cdot f^{\prime })(M_{x}\cdot f)%
\end{array}%
\end{equation*}%
The vertical composition is defined by

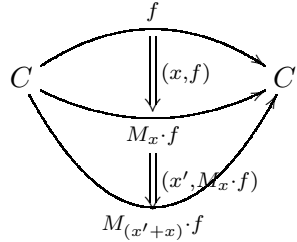
\begin{equation*}
\xymatrix{ C \ar@/^1.6pc/[rrr]_{\quad}^{f}="1" \ar@/_1.2pc/[rrr]_{M_{x}\cdot
f}="2" \ar@/_4pc/[rrr]_{M_{(x'+x)}\cdot f}="9" &&& C
\ar@{}"1";"2"|(.2){\,}="7" \ar@{}"1";"2"|(.8){\,}="8" \ar@{=>}"7"
;"8"^{(x,f)} \ar@{}"2";"9"|(.2){\,}="10" \ar@{}"2";"9"|(.8){\,}="11"
\ar@{=>}"10" ;"11"^{(x',M_{x}\cdot f)} }
\end{equation*}

\begin{equation*}
(x,f)\circ (x^{\prime },M_{x}\cdot f)=(x^{\prime }+x,f)
\end{equation*}%
for $(x,f),(x^{\prime },M_{x}\cdot f)\in C\rtimes \mathcal{M}(C)$ with $%
t(x,f)=s(x^{\prime },M_{x}\cdot f)=M_{x}\cdot f.$

It remains to satisfy the interchange law, i.e.

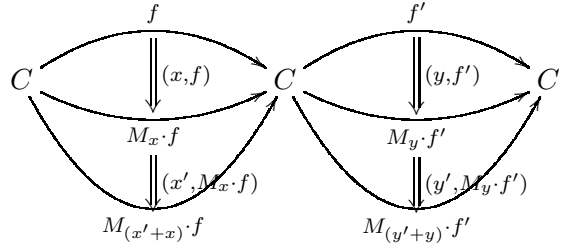
\begin{equation*}
\xymatrix{ C \ar@/^1.6pc/[rrr]_{\quad}^{f}="1" \ar@/_1.2pc/[rrr]_{M_{x}\cdot
f}="2" \ar@/_4pc/[rrr]_{M_{(x'+x)}\cdot f}="9" &&& C
\ar@{}"1";"2"|(.2){\,}="7" \ar@{}"1";"2"|(.8){\,}="8" \ar@{=>}"7"
;"8"^{(x,f)} \ar@{}"2";"9"|(.2){\,}="10" \ar@{}"2";"9"|(.8){\,}="11"
\ar@{=>}"10" ;"11"^{(x',M_{x}\cdot f)} \ar@/^1.6pc/[rrr]_{\quad}^{f'}="12"
\ar@/_1.2pc/[rrr]_{M_{y}\cdot f'}="13"\ar@/_4pc/[rrr]_{M_{(y'+y)}\cdot
f'}="14" \ar@{}"12";"13"|(.2){\,}="15" \ar@{}"12";"13"|(.8){\,}="16"
\ar@{=>}"15" ;"16"^{(y,f')} \ar@{}"13";"14"|(.2){\,}="17"
\ar@{}"13";"14"|(.8){\,}="18" \ar@{=>}"17" ;"18"^{(y',M_{y}\cdot f')} &&& C }
\end{equation*}

\begin{equation*}
\begin{array}{lll}
\lbrack (x,f)\circ (x^{\prime },M_{x}\cdot f)]\bullet \lbrack (y,f^{\prime
})\circ (y^{\prime },M_{y}\cdot f^{\prime })] & = & [(x,f)\bullet
(y,f^{\prime })] \\
&  & \circ \lbrack (x^{\prime },M_{x}\cdot f)\bullet (y^{\prime },M_{y}\cdot
f^{\prime }) ].%
\end{array}%
\end{equation*}

Evaluating the two sides separately, we get
\begin{equation*}
\begin{array}{lll}
\mathbf{LHS} & = & (x^{\prime }+x,f)\bullet (y^{\prime }+y,f^{\prime }) \\
& = & (f^{\prime }(x^{\prime }+x)+f(y^{\prime }+y)+(x^{\prime }+x)(y^{\prime
}+y),f^{\prime }f) \\
& = & (f^{\prime }(x^{\prime })+f^{\prime }(x)+f(y^{\prime })+f(y)+x^{\prime
}y^{\prime }+x^{\prime }y+xy^{\prime }+xy,f^{\prime }f)%
\end{array}%
\end{equation*}%
and

\begin{equation*}
\begin{array}{lll}
\mathbf{RHS} & = & (f^{\prime }(x)+f(y)+xy,f^{\prime }f)\circ ((M_{y}\cdot
f^{\prime })(x^{\prime }) \\
&  & +(M_{x}\cdot f)(y^{\prime })+x^{\prime }y^{\prime },(M_{y}\cdot
f^{\prime })(M_{x}\cdot f)) \\
& = & (f^{\prime }(x)+f(y)+xy+(M_{y}\cdot f^{\prime })(x^{\prime
})+(M_{x}\cdot f)(y^{\prime })+x^{\prime }y^{\prime },f^{\prime }f) \\
& = & (f^{\prime }(x)+f(y)+xy+yx^{\prime }+f^{\prime }(x^{\prime
})+xy^{\prime }+f(y^{\prime })+x^{\prime }y^{\prime },f^{\prime }f)%
\end{array}%
\end{equation*}

LHS and RHS are equal, thus interchange law is satisfied. Therefore we get a
2-algebra consists of the $R$-algebra $C$ as single object and the $R$%
-algebra $A_0$ of 1-morphisms and the $R$-algebra $A_1$ of 2-morphisms.

\section{Crossed modules and 2-algebras}

Crossed modules have been used widely and in various
contexts since their definition by Whitehead \cite{WHITEHEAD} in his investigations of the algebraic structure of
relative homotopy groups. We recalled the definition of crossed modules of
commutative algebras given by Porter \cite{PORTER}.

Let $R$ be a $k$-algebra with identity. A pre-crossed module of commutative
algebras is an $R$-algebra $C$ together with a commutative action of $R$ on $%
C$ and a morphism%
\begin{equation*}
\partial :C\longrightarrow R
\end{equation*}%
such that for all $c\in C$, $r\in R$

\begin{equation*}
\text{CM1) }\partial (r\blacktriangleright c)=r\partial c.
\end{equation*}
This is a crossed $R$-module if in addition for all $c,c^{\prime }\in C$%
\begin{equation*}
\text{CM2) }\partial c\blacktriangleright c^{\prime }=cc^{\prime }.
\end{equation*}%
The last condition is called the Peiffer identity. We denote such a crossed
module by $(C,R,\partial ).$

A morphism of crossed modules from $(C,R,\partial )$ to $(C^{\prime
},R^{\prime },\partial ^{\prime })$ is a pair of $k$-algebra morphisms $\phi
:C\longrightarrow C^{\prime },\psi :R\longrightarrow R^{\prime }$ such that
\begin{equation*}
\partial ^{\prime }\phi =\psi \partial \qquad \text{and}\qquad \phi
(r\blacktriangleright c)=\psi (r)\blacktriangleright \phi (c).
\end{equation*}%
Thus we get a category $\mathbf{XMod}_{k}$ of crossed modules (for fixed $k$%
).

\textbf{Examples of Crossed Modules}

\textbf{1. }Let $I$ be an ideal in $R$. Then $inc:I\longrightarrow R$ is a crossed module. Conversely, if
$\partial :C\longrightarrow R$ is a crossed module then the
Peiffer identity implies that $\partial C$ is an ideal in $R$.

\textbf{2. }Given any $R$-module $M$, the zero morphism $0:M\rightarrow R$ is a crossed module.
Conversely: If $(C,R,\partial )$ is a crossed module, $\partial (C)$ acts
trivially on $\ker \partial ,$ hence $\ker \partial $ has a natural $%
R/\partial (C)$-module structure.

As these two examples suggest, general crossed modules lie between the two
extremes of ideal and modules. Both aspects are important.

\textbf{3.} Let be $\mathcal{M}(C)$ multiplication algebra. Then $\left( C,%
\mathcal{M}\left( C\right) ,\mu \right) $ is multiplication crossed module. $%
\mu :C\rightarrow \mathcal{M}\left( C\right) $ is defined by $\mu \left(
r\right) =\delta _{r}$ with $\delta _{r}\left( r^{\prime }\right)
=rr^{\prime }$ for all $r,r^{\prime }\in C,$ where $\delta $ is multiplier $%
\delta :C\rightarrow C$ such that for all $r,r^{\prime }\in C,$\ $\delta
\left( rr^{\prime }\right) =\delta \left( r\right) r^{\prime }$. Also $%
\mathcal{M}\left( C\right) $ acts on $C$ by $\delta \blacktriangleright
r=\delta \left( r\right) .$(See \cite{ARVASI-EGE} for details).

In \cite{PORTER} Porter states that there is an
equivalence of categories between the category of internal categories in the category of $k$-algebras
 and the category of crossed modules of commutative $k$-algebras. In the following
theorem, we will give a categorical presentation of this equivalence.

\begin{thm}
The category of crossed modules $\mathbf{\ XMod}_{k}$ is equivalent to that
of $2$-algebras, \textbf{2Alg}.
\end{thm}

\begin{pf}
Let $A=(A_{0},A_{1},s,t,e,\circ ,\bullet )$ be a 2-algebra consisting of a
single object say $\ast $ and an algebra $A_{0}$ of 1-morphisms and an
algebra $A_{1}$ of 2-morphisms. For $x,y\in A_{0}$ and $f:x\rightarrow y\in
A_{1}$, we get the following diagram

\begin{equation*}
\xymatrix{ \ast \ar@/^1.5pc/[rr]_{\quad}^{x}="1" \ar@/_1.5pc/[rr]_{y}="2" &&
\ast \ar@{}"1";"2"|(.2){\,}="7" \ar@{}"1";"2"|(.8){\,}="8" \ar@{=>}"7"
;"8"^{f} }
\end{equation*}

We define $s,t$ morphisms $s:A_{1}\longrightarrow
A_{0},s(f)=x,t:A_{1}\longrightarrow A_{0},t(f)=y$ and $e$ morphism $%
e:A_{0}\longrightarrow A_{1}$ for $x\in A_{0},$ $e(x):x\longrightarrow x\in
A_{1}$.

The $s,t$ and $e$ morphisms are algebra morphisms and we have
\begin{eqnarray*}
se(x) &=&s(e(x))=x=Id_{A_{0}}(x) \\
te(x) &=&t(e(x))=x=Id_{A_{0}}(x)
\end{eqnarray*}%
We define
\begin{equation*}
\text{Ker }s=K=\{q\in A_{1}\mid s(q)=Id_{A_{0}}\}\subseteq A_{1}
\end{equation*}%
and $\partial =t\mid _{K}$ algebra homomorphism by $\partial
:K\longrightarrow A_{0},\partial (q)=t(q)$. \ We have semidirect product Ker
$s\rtimes A_{0}=\{(q,x)\mid q\in $Ker$s,x\in A_{0}\}$ with multiplication $%
(q,x)\bullet (q^{\prime },x^{\prime })=(x\blacktriangleright q^{\prime
}+x^{\prime }\blacktriangleright q+q^{\prime }\bullet q,x\bullet x^{\prime
}) $ where action of $A_{0}$ on Ker$s$ is defined by $x\blacktriangleright
q=e(x)\bullet h$. For each $f\in A_{1}$, we can write $f=q+e(x)$ where $%
q=f-es(f)\in $Ker$s$ and $x=s(f).$ Suppose $f^{\prime }=q^{\prime
}+e(x^{\prime }).$ Then
\begin{equation*}
\begin{array}{lll}
f\bullet f^{\prime } & = & (q+e(x))\bullet (q^{\prime }+e(x^{\prime })) \\
& = & q\bullet q^{\prime }+q\bullet e(x^{\prime })+e(x)\bullet q^{\prime
}+e(x)\bullet e(x^{\prime }) \\
& = & e(x^{\prime })\bullet q+e(x)\bullet q^{\prime }+q\bullet q^{\prime
}+e(x\bullet x^{\prime }) \\
& = & x^{\prime }\blacktriangleright q+x\blacktriangleright q^{\prime
}+q\bullet q^{\prime }+e(x\bullet x^{\prime }).%
\end{array}%
\end{equation*}%
There is a map%
\begin{equation*}
\begin{array}{cccc}
\phi : & A_{1} & \longrightarrow & \text{Ker}s\rtimes A_{0} \\
& q+e(x) & \longmapsto & \phi (q+e(x))=(q,x).%
\end{array}%
\end{equation*}%
Now
\begin{equation*}
\begin{array}{lll}
\phi (f\bullet f^{\prime }) & = & \phi (x^{\prime }\blacktriangleright
q+x\blacktriangleright q^{\prime }+q\bullet q^{\prime }+e(x\bullet x^{\prime
})) \\
& = & (x^{\prime }\blacktriangleright q+x\blacktriangleright q^{\prime
}+q\bullet q^{\prime },x\bullet x^{\prime }) \\
& = & (q,x)\bullet (q^{\prime },x^{\prime }) \\
& = & \phi (f)\bullet \phi (f^{\prime })%
\end{array}%
\end{equation*}%
so $\phi $ is a homomorphism. Also, there is an obvious inverse%
\begin{equation*}
\begin{array}{cccc}
\phi ^{-1}: & \text{Ker}s\rtimes A_{0} & \longrightarrow & A_{1} \\
& (q,x) & \longmapsto & \phi ^{-1}(q,x)=q+e(x)%
\end{array}%
\end{equation*}%
which is also a homomorphism. Hence $\phi $ is an isomorphism and we have
established that Ker $s\rtimes A_{0}\simeq A_{1}$. Since $A$ is a $2$%
-algebra and Ker $s\rtimes A_{0}\simeq A_{1}$, we can define algebra
morphisms%
\begin{equation*}
\begin{array}{rcl}
s: & \text{Ker}s\rtimes A_{0} & \longrightarrow A_{0} \\
& (q,x) & \longmapsto s(q,x)=x%
\end{array}%
\ \ \ \ \ \ \
\begin{array}{rcl}
t: & \text{Ker}s\rtimes A_{0} & \longrightarrow A_{0} \\
& (q,x) & \longmapsto t(q,x)=\partial (q)+x%
\end{array}%
\end{equation*}%
and%
\begin{equation*}
\begin{array}{rcl}
e: & A_{0} & \longrightarrow \text{Ker}s\rtimes A_{0} \\
& x & \longmapsto e(x)=(0,x)%
\end{array}%
\end{equation*}%
and for $t(q,x)=s(q^{\prime },\partial (q)+x)=\partial (q)+x$ we define
\begin{equation*}
\begin{array}{ccc}
\circ : & \text{Ker}s\rtimes A_{0~t}\times _{~s}\text{Ker}s\rtimes A_{0} &
\longrightarrow \text{Ker}s\rtimes A_{0} \\
& \left( (q,x),(q^{\prime },\partial (q)+x)\right) & \longmapsto (q^{\prime
}+q,x)%
\end{array}%
\end{equation*}

\begin{equation*}
\xymatrix{ \ast \ar@/^1.6pc/[rrr]_{\quad}^{x}="1"
\ar@/_1.2pc/[rrr]_{\partial(q)+x}="2" \ar@/_4pc/[rrr]_{\partial(q'+q)+x}="9"
&&& \ast \ar@{}"1";"2"|(.2){\,}="7" \ar@{}"1";"2"|(.8){\,}="8" \ar@{=>}"7"
;"8"^{(q,x)} \ar@{}"2";"9"|(.2){\,}="10" \ar@{}"2";"9"|(.8){\,}="11"
\ar@{=>}"10" ;"11"^{(q',\partial(q)+x)} & =} \xymatrix{ \qquad \ast
\ar@/^1.5pc/[rrrr]_{\quad}^{x}="1" \ar@/_1.5pc/[rrrr]_{\partial(q'+q)+x}="2"
&&&& \ast \ar@{}"1";"2"|(.2){\,}="7" \ar@{}"1";"2"|(.8){\,}="8" \ar@{=>}"7"
;"8"^{(q'+q,x)} }
\end{equation*}

which is vertical composition;%
\begin{equation*}
(q,x)\circ (q^{\prime },\partial (q)+x)=(q^{\prime }+q,x).
\end{equation*}%
For\textbf{\ }$(q,x),(p,y)\in $Ker$s\rtimes A_{0}$, horizontal composition
is defined by

\begin{equation*}
\xymatrix{ \ast \ar@/^2pc/[rr]_{\quad}^{x}="1"
\ar@/_2pc/[rr]_{\partial(q)+x}="2" && \ast\ar@{}"1";"2"|(.2){\,}="7"
\ar@{}"1";"2"|(.8){\,}="8" \ar@{=>}"7" ;"8"^{(q,x)}
\ar@/^2pc/[rr]_{\quad}^{y}="9" \ar@/_2pc/[rr]_{\partial(p)+y}="10"
\ar@{}"9";"10"|(.2){\,}="11" \ar@{}"9";"10"|(.8){\,}="12" \ar@{=>}"11"
;"12"^{(p,y)} && \ast & = } \xymatrix{ \qquad \ast
\ar@/^1.5pc/[rrrrr]_{\quad}^{x\bullet y}="1"
\ar@/_1.5pc/[rrrrr]_{(\partial(q)+x)\bullet (\partial(p)+y)}="2" &&&&&\ast
\ar@{}"1";"2"|(.2){\,}="7" \ar@{}"1";"2"|(.8){\,}="8" \ar@{=>}"7"
;"8"^{(x\blacktriangleright p +y\blacktriangleright q+p\bullet q,x\bullet
y)} }
\end{equation*}

\begin{equation*}
\begin{array}{ccc}
(q,x)\bullet (p,y) & = & (x\blacktriangleright p+y\blacktriangleright
q+p\bullet q,x\bullet y) \\
& = & (e(x)\bullet p+e(y)\bullet q+p\bullet q,x\bullet y).%
\end{array}%
\end{equation*}

Thus we have

CM1)

$%
\begin{array}{lll}
\partial (x\blacktriangleright q) & = & \partial (e(x)\bullet q) \\
& = & \partial (e(x))\bullet \partial (q) \\
& = & \left( te\right) (x)\bullet \partial (q) \\
& = & x\bullet \partial (q).%
\end{array}%
$

Also by interchange law we have
\begin{equation*}
\begin{array}{ccc}
\left[ (q,x)\bullet (p,y)\right] \circ \left[ (q^{\prime },\partial
(q)+x)\bullet (p^{\prime },\partial (p)+y)\right] & = & \left[ (q,x)\circ
(q^{\prime },\partial (q)+x)\right] \\
&  & \bullet \left[ (p,y)\circ (p^{\prime },\partial (p)+y)\right] .%
\end{array}%
\end{equation*}

Therefore, evaluating the two sides of this equation gives:
\begin{equation*}
\begin{array}{lll}
LHS & = & (x\blacktriangleright p+y\blacktriangleright q+p\bullet q,x\bullet
y) \\
&  & \circ (\left( \partial (q)+x\right) \blacktriangleright p^{\prime
}+\left( \partial (p)+y\right) \blacktriangleright q^{\prime }+p^{\prime
}\bullet q^{\prime },\left( \partial (q)+x\right)\bullet \left( \partial
(p)+y\right) ) \\
& = & (\left( \partial (q)+x\right) \blacktriangleright p^{\prime }+\left(
\partial (p)+y\right) \blacktriangleright q^{\prime }+p^{\prime }\bullet
q^{\prime }+x\blacktriangleright p+y\blacktriangleright q+p\bullet
q,x\bullet y) \\
& = & (\partial (q)\blacktriangleright p^{\prime }+e(x)\bullet p^{\prime
}+\partial (p)\blacktriangleright q^{\prime } \\
&  & +e(y)\bullet q^{\prime }+p^{\prime }\bullet q^{\prime }+e(x)\bullet
p+e(y)\bullet q+p\bullet q,x\bullet y) \\
RHS & = & (q^{\prime }+q,x)\bullet (p^{\prime }+p,y) \\
& = & (x\blacktriangleright \left( p^{\prime }+p\right)
+y\blacktriangleright (q^{\prime }+q)+\left( p^{\prime }+p\right)\bullet
(q^{\prime }+q),x\bullet y) \\
& = & \left( e(x)\bullet p^{\prime }+e(x)\bullet p+e(y)\bullet q^{\prime
}+e(y)\bullet q+p^{\prime }\bullet q^{\prime }+p^{\prime }\bullet q+p\bullet
q^{\prime }+p\bullet q,x\bullet y\right) .%
\end{array}%
\end{equation*}

Since the two sides are equal, we know that their first components must be
equal, so we have
\begin{equation*}
\partial (q)\blacktriangleright p^{\prime }+\partial (p)\blacktriangleright
q^{\prime }=q\bullet p^{\prime }+p\bullet q^{\prime }
\end{equation*}%
and%
\begin{equation*}
\begin{array}{lll}
q\bullet p^{\prime }+p\bullet q^{\prime } & = & \partial
(q)\blacktriangleright p^{\prime }+\partial (p)\blacktriangleright q^{\prime
} \\
& = & \partial (q+p)\blacktriangleright (p^{\prime }+q^{\prime })-\partial
(q)\blacktriangleright q^{\prime }-\partial (p)\blacktriangleright p^{\prime
} \\
& = & \partial (q+p)\blacktriangleright (p^{\prime }+q^{\prime })-\left(
q\bullet q^{\prime }+p\bullet p^{\prime }\right) ,%
\end{array}%
\end{equation*}%
thus
\begin{equation*}
\begin{array}{lll}
\partial (q+p)\blacktriangleright (p^{\prime }+q^{\prime }) & = & q\bullet
p^{\prime }+p\bullet q^{\prime }+\left( q\bullet q^{\prime }+p\bullet
p^{\prime }\right)  \\
& = & \left( q+p\right) \bullet \left( q^{\prime }+p^{\prime }\right)
\end{array}%
\end{equation*}%
and writing $(q+p)=l,\left( q^{\prime }+p^{\prime }\right) =l^{\prime }\in
Kers,$ we get
\begin{equation*}
\partial \left( l\right) \blacktriangleright l^{\prime }=l\bullet l^{\prime }
\end{equation*}%
which is the Peiffer identity as required. Hence $\left( Kers,A_{0},\partial
\right) $ is a crossed module.

Let $\mathcal{A}=(A_{0},A_{1},s,t,e,\circ ,\bullet )$ and $\mathcal{A}%
^{^{\prime }}=(A_{0}^{^{\prime }},A_{1}^{^{\prime }},s^{^{\prime
}},t^{^{\prime }},e^{^{\prime }},\circ ^{^{\prime }},\bullet ^{^{\prime }})$
be 2-algebras and $F=(F_{0},F_{1}):$ $\mathcal{A\longrightarrow A}^{^{\prime
}}$ be a 2-algebra morphism. Then $F_{0}:A_{0}\longrightarrow
A_{0}^{^{\prime }}$ and $F_{1}:A_{1}\longrightarrow A_{1}^{^{\prime }}$ are
the $k$-algebra morphisms. We define $f_{1}=F_{1}|_{Kers}:Kers%
\longrightarrow Kers^{^{\prime }}$ and $f_{0}=F_{0}:A_{0}\longrightarrow
A_{0}^{^{\prime }}.$ For all $a\in Kers$ and $x\in A_{0}$,%
\begin{equation*}
\begin{array}{lll}
f_{0}\partial (a) & = & F_{0}t(a) \\
& = & t^{^{\prime }}F_{1}(a) \\
& = & \partial ^{^{\prime }}f_{1}(a)%
\end{array}%
\end{equation*}%
and%
\begin{equation*}
\begin{array}{lll}
f_{1}(x\blacktriangleright a) & = & F_{1}(e(x)a) \\
& = & F_{1}(e(x))F_{1}(a) \\
& = & e^{^{\prime }}F_{0}(x)F_{1}(a) \\
& = & e^{^{\prime }}f_{0}(x)f_{1}(a) \\
& = & f_{0}(x)\blacktriangleright f_{1}(a).%
\end{array}%
\end{equation*}%
Thus $(f_{1},f_{0})$ map is a crossed module morphism $(Kers,A_{0},\partial
)\longrightarrow (Kers^{^{\prime }},A_{0}^{^{\prime }},\partial ^{^{\prime
}}).$ So we have a functor%
\begin{equation*}
\Gamma :\mathbf{2Alg}\longrightarrow \mathbf{XMod}_{k}.
\end{equation*}

Conversely, let $(G,C,\partial )$ be a crossed module of algebras. Therefore
there is an algebra morphism $\partial :G\rightarrow C$ and an action of $C$
on $G$ such that

CM1) $\partial (x\blacktriangleright g)=x\partial (g),$

CM2) $\partial (g)\blacktriangleright g^{\prime }=gg^{\prime }.$

Since $C$ acts on $G$, we can form the semidirect product $G\rtimes C$ as
defined by
\begin{equation*}
G\rtimes C=\{\left( g,c\right) \mid g\in G,c\in C\}
\end{equation*}%
with multiplication%
\begin{equation*}
\left( g,c\right) \left( g^{\prime },c^{\prime }\right) =\left(
c\blacktriangleright g^{\prime }+c^{\prime }\blacktriangleright g+g^{\prime
}g,cc^{\prime }\right)
\end{equation*}%
and define maps $s,t:G\rtimes C\rightarrow C$ and $e:C\rightarrow G\rtimes C$
by $s(g,c)=c,~t(g,c)=\partial (g)+c$ and $e(c)=(0,c).$ These maps are
clearly algebra morphisms.

\begin{equation*}
\xymatrix{ \ast \ar@/^1.6pc/[rrrr]_{\quad}^{c}="1"
\ar@/_1.2pc/[rrrr]_{\partial(g)+c}="2"
\ar@/_4pc/[rrrr]_{\partial(g+g')+c}="9" &&&& \ast \ar@{}"1";"2"|(.2){\,}="7"
\ar@{}"1";"2"|(.8){\,}="8" \ar@{=>}"7" ;"8"^{(g,c)}
\ar@{}"2";"9"|(.2){\,}="10" \ar@{}"2";"9"|(.8){\,}="11" \ar@{=>}"10"
;"11"^{(g',\partial(g)+c)} }
\end{equation*}
For $t(g,c)=s(g^{\prime },\partial (g)+c)=\partial (g)+c,$ we define
composition

\begin{equation*}
\begin{array}{ccl}
\circ : & \left( G\rtimes C\right) _{t}\times _{~s}\left( G\rtimes C\right)
& \longrightarrow \left( G\rtimes C\right) \\
& \left( (g,c),(g^{\prime },\partial (g)+c)\right) & \longmapsto
(g+g^{\prime },c),%
\end{array}%
\end{equation*}%
for $(g,c),(h,d)\in G\rtimes C$ and $(g,c),(g^{\prime },\partial (g)+c)\in
G\rtimes C,$ following equations give horizontal and vertical composition
respectively.
\begin{equation*}
(g,c)\bullet (h,d)=(c\blacktriangleright h+d\blacktriangleright g+gh,cd)
\end{equation*}%
\begin{equation*}
(g,c)\circ (g^{\prime },\partial (g)+c)=(g+g^{\prime },c)
\end{equation*}%
Finally, since it must be that $\circ $ is an algebra morphism and by the
crossed module conditions, interchange law is satisfied. Therefore we have
constructed a $2$-algebra $\mathcal{A}=(C,G\rtimes C,s,t,e,\circ ,\bullet )$
consists of the single object say $\ast $ and the $k$-algebra $C$ of
1-morphisms and the $k$-algebra $G\rtimes C$ of 2-morphisms. Let $%
(G,C,\partial )$ and $(G^{^{\prime }},C^{^{\prime }},\partial ^{^{\prime }})$
be crossed modules and $f=(f_{1},f_{0}):(G,C,\partial )$ $\longrightarrow
(G^{^{\prime }},C^{^{\prime }},\partial ^{^{\prime }})$ be a crossed module
morphism. We define
\begin{equation*}
\begin{array}{cccc}
F_{1}: & G\rtimes C & \longrightarrow & G^{^{\prime }}\rtimes
C^{^{\prime }} \\
& (g,c) & \longmapsto & F_{1}(g,c)=(f_{1}(g),f_{0}(c))%
\end{array}%
\end{equation*}%
and
\begin{equation*}
\begin{array}{cccc}
F_{0}: & C & \longrightarrow & C^{^{\prime }} \\
& c & \longmapsto & F_{0}(c)=f_{0}(c).%
\end{array}%
\end{equation*}%
Then
\begin{equation*}
\begin{array}{lll}
s^{^{\prime }}F_{1}(g,c) & = & s^{^{\prime }}(f_{1}(g),f_{0}(c)) \\
& = & f_{0}(c) \\
& = & F_{0}(c) \\
& = & F_{0}s(g,c),%
\end{array}%
\end{equation*}%
\begin{equation*}
\begin{array}{lll}
t^{^{\prime }}F_{1}(g,c) & = & t^{^{\prime }}(f_{1}(g),f_{0}(c)) \\
& = & \partial ^{^{\prime }}f_{1}(g)+f_{0}(c) \\
& = & f_{0}\partial (g)+f_{0}(c) \\
& = & F_{0}(\partial (g)+c) \\
& = & F_{0}t(g,c),%
\end{array}%
\end{equation*}

\begin{equation*}
\begin{array}{lll}
e^{^{\prime }}F_{0}(c) & = & (0,f_{0}(c)) \\
& = & F_{1}(0,c) \\
& = & F_{1}e(c),%
\end{array}%
\end{equation*}%
\begin{equation*}
\begin{array}{lll}
F_{1}((g,c)\circ (g^{^{\prime }},c^{^{\prime }})) & = & F_{1}(g+g^{^{\prime
}},c) \\
& = & (f_{1}(g+g^{^{\prime }}),f_{0}(c)) \\
& = & (f_{1}(g)+f_{1}(g^{^{\prime }}),f_{0}(c)) \\
& = & (f_{1}(g),f_{0}(c))\circ (f_{1}(g^{^{\prime }}),f_{0}(c^{^{\prime }}))
\\
& = & F_{1}(g,c)\circ F_{1}(g^{^{\prime }},c^{^{\prime }}),%
\end{array}%
\end{equation*}

\begin{equation*}
\begin{array}{lll}
F_{1}((g,c)\bullet (h,d)) & = & F_{1}(c\blacktriangleright
h+d\blacktriangleright g+gh,cd) \\
& = & (f_{1}(c\blacktriangleright h)+f_{1}(d\blacktriangleright
g)+f_{1}(gh),f_{0}(cd)) \\
& = & (f_{0}(c)\blacktriangleright f_{1}(h)+f_{0}(d)\blacktriangleright
f_{1}(g)+f_{1}(g)f_{1}(h),f_{0}(c)f_{0}(d)) \\
& = & (f_{1}(g),f_{0}(c))\bullet (f_{1}(h),f_{0}(d)) \\
& = & F_{1}(g,c)\bullet F_{1}(h,d)%
\end{array}%
\end{equation*}
for all $(g,c)\in G\rtimes C$ and $c\in C$. Therefore $%
F=(F_{1},F_{0})$ is a 2-algebra morphism from $%
(C,G\rtimes C,s,t,e,\circ ,\bullet )$ to $(C^{^{\prime }},G^{^{\prime
}}\rtimes C^{^{\prime }},s^{^{\prime }},t^{^{\prime }},e^{^{\prime }},\circ
^{^{\prime }},\bullet ^{^{\prime }}).$ Thus we get a functor
\begin{equation*}
\Psi :\mathbf{XMod}_{k}\longrightarrow \mathbf{2Alg}.
\end{equation*}
\end{pf}

\subsection{ Homotopies of Crossed modules and 2-algebras}

The notion of homotopy for morphisms of crossed modules over commutative
algebras is given in \cite{kim}. In this section, we explain the relation
between homotopies for crossed modules over commutative algebras and
homotopies for 2-algebras. The formulae given below are playing important
role in our study.

\begin{defn}
\cite{kim} Let $\mathcal{A}=(E,R,\partial )$ and $\mathcal{A}^{^{\prime
}}=(E^{^{\prime }},R^{^{\prime }},\partial ^{^{\prime }})$ be crossed
modules and $f_{0}:R\longrightarrow R^{^{\prime }}$ be an algebra morphism.
An $f_{0}$-derivation $s:R\longrightarrow E^{^{\prime }}$ is a $k$-linear
map satisfying for all $r,r^{^{\prime }}\in R,$%
\begin{equation*}
s(rr^{^{\prime }})=f_{0}(r)\blacktriangleright s(r^{^{\prime
}})+f_{0}(r^{^{\prime }})\blacktriangleright s(r)+s(r)s(r^{^{\prime }}).
\end{equation*}%
Let $f=(f_{1},f_{0})$ be a crossed module morphism $\mathcal{A}%
\longrightarrow \mathcal{A}^{^{\prime }}$ and $s$ be an $f_{0}$-derivation. If $g=(g_{1},g_{2})$ is defined as (where $e\in E$ and $r\in R$)%
\begin{equation*}
\begin{array}{ccc}
g_{0}(r) & = & f_{0}(r)+(\partial ^{^{\prime }} s)(r) \\
g_{1}(e) & = & f_{1}(e)+(s \partial )(e),%
\end{array}%
\end{equation*}%
then $g$ is also crossed module morphism $\mathcal{A}\longrightarrow
\mathcal{A}^{^{\prime }}.$ In such a case we write $f\overset{(f_{0},s)}{%
\longrightarrow }g$, and say that $(f_{0},s)$ is a homotopy connecting $f$
to $g$.
\end{defn}

If $(f_{0},s)$ and $(g_{0},s^{^{\prime }})$ are homotopies connecting $f$ to
$g$ and $g$ to $u$ respectively, then $(f_{0},s+s^{^{\prime }})$ is a
homotopy connecting $f$ to $u$, where $s+s^{^{\prime }}:R\longrightarrow
E^{^{\prime }}$ is an $f_{0}$-derivation defined by $(s+s^{^{\prime
}})(r)=s(r)+s^{^{\prime }}(r)$.

The notion of homotopy for $2$-algebras is essentially a special case of $2$%
-natural transformation due to Gray in \cite{Gray}.

\begin{defn}
Let $\mathbf{A=(}A_{0},A_{1},s,t,e,\circ ,\bullet \mathbf{)}$ and $\mathbf{A}%
^{^{\prime }}\mathbf{=(}A_{0}^{^{\prime }},A_{1}^{^{\prime }},s^{^{\prime
}},t^{^{\prime }},e^{^{\prime }},\circ ^{^{\prime }},{\bullet }^{^{\prime }}%
\mathbf{)}$ be 2-algebras and let $F=(F_{1},F_{0})$ and $G=(G_{1},G_{0})$ be
2-algebra morphisms $\mathbf{A}\longrightarrow \mathbf{A}^{^{\prime }}.$ A
k-algebra morphism $\delta :A_{0}\longrightarrow A_{1}^{^{\prime }}$ satisfying
the following conditions is called a homotopy connecting $F$ to $G:$

1) $s^{^{\prime }} \delta =F_{0}$

2) $t^{^{\prime }} \delta =G_{0}$

3) $F_{1}\circ ^{^{\prime }}\delta t=\delta s\circ ^{^{\prime }}G_{1}.$ In such a
case we write $F\overset{\delta }{\longrightarrow }G$.
\end{defn}

\begin{thm}
Let $\mathcal{A}=(A_{0},A_{1},s,t,e,\circ ,\bullet )$ , $\mathcal{A}%
^{^{\prime }}=(A_{0}^{^{\prime }},A_{1}^{^{\prime }},s^{^{\prime
}},t^{^{\prime }},e^{^{\prime }},\circ ^{^{\prime }},\bullet ^{^{\prime }})$
be 2-algebras, $F=(F_{1},F_{0})$, $G=(G_{1},G_{0})$ and $U=(U_{1},U_{0})$ be
2-algebra morphisms $\mathcal{A\longrightarrow A}^{^{\prime }}$and $\delta $ be
a homotopy connecting $F$ to $G,$ $\delta ^{^{\prime }}$ be a homotopy
connecting $G$ to $U$. Then the map $\delta \ast \delta ^{^{\prime
}}:A_{0}\longrightarrow A_{1}$ defined by $(\delta \ast \delta ^{^{\prime
}})(x)=\delta (x)+\delta ^{^{\prime }}(x)-e^{^{\prime }}(t^{^{\prime }}\delta )(x)$
is a homotopy connecting $F$ to $U$.
\end{thm}

\begin{pf}
We first show that $\delta \ast \delta ^{^{\prime }}$ is an algebra morphism.
Since $\delta $ and $\delta ^{^{\prime }}$ are algebra morphisms, $\delta (x\bullet
x^{^{\prime }})=\delta (x)\bullet ^{^{\prime }}\delta (x^{^{\prime }})$ and $\delta
^{^{\prime }}(x\bullet x^{^{\prime }})=\delta ^{^{\prime }}(x)\bullet
^{^{\prime }}\delta ^{^{\prime }}(x^{^{\prime }})$ for all $x,x^{^{\prime }}\in
A_{0}.$ Then we get
\begin{equation*}
\begin{array}{lll}
(\delta \ast \delta ^{^{\prime }})(x\bullet x^{^{\prime }}) & = & \delta (x\bullet
x^{^{\prime }})+\delta ^{^{\prime }}(x\bullet x^{^{\prime }})-e^{^{\prime
}}(t^{^{\prime }}\delta )(x\bullet x^{^{\prime }}) \\
& = & \delta (x\bullet x^{^{\prime }})+\delta ^{^{\prime }}(x\bullet x^{^{\prime
}})-e^{^{\prime }}(G_{0})(x\bullet x^{^{\prime }}) \\
& = & \delta (x\bullet x^{^{\prime }})\circ ^{^{\prime }}\delta ^{^{\prime
}}(x\bullet x^{^{\prime }})\text{ \ \ \ \ \ \ \ \ \ \ \ \ \ (Proposition \ref{2.2})} \\
& = & (\delta (x)\bullet ^{^{\prime }}\delta (x^{^{\prime }}))\circ ^{^{\prime
}}(\delta ^{^{\prime }}(x)\bullet ^{^{\prime }}\delta ^{^{\prime }}(x^{^{\prime
}})) \\
& = & (\delta (x)\circ ^{^{\prime }}\delta ^{^{\prime }}(x))\bullet ^{^{\prime
}}(\delta (x^{^{\prime }})\circ ^{^{\prime }}\delta ^{^{\prime }}(x^{^{\prime }}))%
\text{ \ \ \ \ \ \ \ \ (interchange law)} \\
& = & (\delta (x)+\delta ^{^{\prime }}(x)-e^{^{\prime }}(G_{0})(x))\bullet
^{^{\prime }}(\delta (x^{^{\prime }})+\delta ^{^{\prime }}(x^{^{\prime
}})-e^{^{\prime }}(G_{0})(x^{^{\prime }})) \\
& = & (\delta \ast \delta ^{^{\prime }})(x)\bullet ^{^{\prime }}(\delta \ast \delta
^{^{\prime }})(x^{^{\prime }}).%
\end{array}%
\end{equation*}%
For all $x\in A_{0}$%
\begin{equation*}
\begin{array}{lll}
s^{^{\prime }}(\delta \ast \delta ^{^{\prime }})(x) & = & s^{^{\prime }}(\delta
(x)+\delta ^{^{\prime }}(x)-e^{^{\prime }}G_{0}(x)) \\
& = & s^{^{\prime }}\delta (x)+s^{^{\prime }}\delta ^{^{\prime }}(x)-s^{^{\prime
}}e^{^{\prime }}G_{0}(x) \\
& = & F_{0}(x)+G_{0}(x)-G_{0}(x) \\
& = & F_{0}(x),%
\end{array}%
\end{equation*}

\begin{equation*}
\begin{array}{lll}
t^{^{\prime }}(\delta \ast \delta ^{^{\prime }})(x) & = & t^{^{\prime }}(\delta
(x)+\delta ^{^{\prime }}(x)-e^{^{\prime }}G_{0}(x)) \\
& = & t^{^{\prime }}\delta (x)+t^{^{\prime }}\delta ^{^{\prime }}(x)-t^{^{\prime
}}e^{^{\prime }}G_{0}(x) \\
& = & G_{0}(x)+U_{0}(x)-G_{0}(x) \\
& = & U_{0}(x),%
\end{array}%
\end{equation*}%
and since $F_{1}\circ ^{^{\prime }}\delta t$ $=\delta s\circ ^{^{\prime }}G_{1}$
and $G_{1}\circ ^{^{\prime }}\delta ^{^{\prime }}t$ $=\delta ^{^{\prime }}s\circ
^{^{\prime }}U_{1}$, we get
\begin{equation*}
\begin{array}{ccc}
F_{1}\circ ^{^{\prime }}\delta t\circ ^{^{\prime }}\delta ^{^{\prime }}t & = & \delta
s\circ ^{^{\prime }}G_{1}\circ ^{^{\prime }}\delta ^{^{\prime }}t \\
& = & \delta s\circ ^{^{\prime }}\delta ^{^{\prime }}s\circ ^{^{\prime }}U_{1}.%
\end{array}%
\end{equation*}%
Thus, we get
\begin{equation*}
\begin{array}{ccc}
F_{1}\circ ^{^{\prime }}(\delta \ast \delta ^{^{\prime }})t & = & F_{1}\circ
^{^{\prime }}(\delta t\circ ^{^{\prime }}\delta ^{^{\prime }}t) \\
& = & (\delta s\circ ^{^{\prime }}\delta ^{^{\prime }}s)\circ ^{^{\prime }}U_{1}
\\
& = & (\delta \ast \delta ^{^{\prime }})s\circ ^{^{\prime }}U_{1}.%
\end{array}%
\end{equation*}

Therefore $\delta \ast \delta ^{^{\prime }}:A_{0}\longrightarrow A_{1}$ is a
homotopy connecting $F$ to $U.$
\end{pf}

\begin{thm}
Let $\Gamma :2A\lg \longrightarrow XMod_{k}$ be the functor as mentioned in
Theorem 3.1 and $\delta $ be homotopy connecting $F$ to $G$. Then
\begin{equation*}
\begin{array}{ccccc}
\Gamma (\delta )=h & : & A_{0} & \longrightarrow & Kers^{^{\prime }} \\
&  & x & \longmapsto & h(x)=\delta (x)-e^{^{\prime }}(s^{^{\prime }}\delta )(x)%
\end{array}%
\end{equation*}%
is a homotopy of corresponding crossed module morphisms.
\end{thm}

\begin{pf}
We first show that $h$ is an $f_{0}-$derivation where $f_{0}:A_{0}%
\longrightarrow A_{0}^{^{\prime }}$ defined by $f_{0}(x)=F_{0}(x)$. For $%
x,x^{^{\prime }}\in A_{0},$%
\begin{equation*}
\begin{array}{lll}
f_{0}(x)\blacktriangleright h(x^{^{\prime }}) &  &  \\
+f_{0}(x^{^{\prime }})\blacktriangleright h(x)+h(x)\bullet ^{^{\prime
}}h(x^{^{\prime }}) & = & F_{0}(x)\blacktriangleright (\delta (x^{^{\prime
}})-e^{^{\prime }}(s^{^{\prime }}\delta )(x^{^{\prime }})) \\
&  & +F_{0}(x^{^{\prime }})\blacktriangleright (\delta (x)-e^{^{\prime
}}(s^{^{\prime }}\delta )(x)) \\
&  & +(\delta (x)-e^{^{\prime }}(s^{^{\prime }}\delta )(x))\bullet ^{^{\prime
}}(\delta (x^{^{\prime }})-e^{^{\prime }}(s^{^{\prime }}\delta )(x^{^{\prime }}))
\\
& = & e^{^{\prime }}(F_{0}(x))\bullet ^{^{\prime }}(\delta (x^{^{\prime
}})-e^{^{\prime }}F_{0}(x^{^{\prime }})) \\
&  & +e^{^{\prime }}(F_{0}(x^{^{\prime }}))\bullet ^{^{\prime }}(\delta
(x)-e^{^{\prime }}F_{0}(x))+\delta (x)\bullet ^{^{\prime }}\delta (x^{^{\prime }})
\\
&  & -\delta (x)\bullet ^{^{\prime }}e^{^{\prime }}F_{0}(x^{^{\prime
}})-e^{^{\prime }}F_{0}(x)\bullet ^{^{\prime }}\delta (x^{^{\prime
}})+e^{^{\prime }}F_{0}(x)\bullet ^{^{\prime }}e^{^{\prime
}}F_{0}(x^{^{\prime }}) \\
& = & \delta (x\bullet x^{^{\prime }})-e^{^{\prime }}(s^{^{\prime }}\delta
)(x\bullet x^{^{\prime }}) \\
& = & h(x\bullet x^{^{\prime }}).%
\end{array}%
\end{equation*}%
Therefore $h$ is an $f_{0}-$derivation.

Now we show that
\begin{equation*}
\begin{array}{ccc}
g_{0}(x) & = & f_{0}(x)+\partial ^{^{\prime }}h(x) \\
g_{1}(n) & = & f_{1}(n)+h\partial (n)%
\end{array}%
\end{equation*}%
for $x\in A_{0}$ and $n\in Kers.$%
\begin{equation*}
\begin{array}{lll}
\partial ^{^{\prime }}h(x) & = & \partial ^{^{\prime }}(\delta (x)-e^{^{\prime
}}f_{0}(x)) \\
& = & \partial ^{^{\prime }}(\delta (x))-\partial ^{^{\prime }}(e^{^{\prime
}}f_{0}(x)) \\
& = & (t^{^{\prime }}\delta )(x)-(t^{^{\prime }}e^{^{\prime }})f_{0}(x) \\
& = & g_{0}(x)-f_{0}(x)%
\end{array}%
\end{equation*}%
and we get $g_{0}(x)=f_{0}(x)+\partial ^{^{\prime }}h(x)$.

Since $A_{1}\simeq Kers\rtimes A_{0}$, we take $a=(n,x)$ for $a\in A_{1}$
where $n=a-es(a)\in Kers$ and $x=s(a)\in A_{0}.$ We define $\delta ^{\ast
}:A_{0}\longrightarrow $\ $Kers^{^{\prime }}\rtimes A_{0}^{^{\prime }}$, as $%
\delta ^{\ast }(x)=(\delta (x)-e^{^{\prime }}s^{^{\prime }}(\delta (x)),s^{^{\prime
}}\delta (x))$ and $h^{\ast }:A_{0}\longrightarrow $\ $Kers^{^{\prime }}\rtimes
A_{0}^{^{\prime }}$, as $h^{\ast }(x)=(h(x),F_{0}(x))$. Therefore
\begin{equation*}
\xymatrix@R=60pt@C=60pt{
	*+[l] {A_1 \cong Ker(s) \rtimes A_0}
	\ar[r]^-{s} \ar@<-1.5ex>[r]_-{t}
	\ar@/^0.5pc/[d]^{(F_{1}, F_0)}
	\ar@/_0.5pc/[d]_{(G_{1}, G_0)}
	& A_{0}
	\ar@{}[dl]^(.15){}="a"^(.85){}="b" \ar@[red] "a";"b"^-{\delta^{\ast}}
	\ar@/^0.5pc/[d]^{F_{0}}
	\ar@/_0.5pc/[d]_{G_{0}}
	\ar@/_1.5pc/[l]_{e}
	\\
	*+[l]{A'_1 \cong Ker(s') \rtimes A'_0}	\ar@<1.5ex>[r]^-{s'} \ar[r]_-{t'}
	& A'_{0} \ar@/^1.5pc/[l]^{e'} }
\end{equation*}
for $(F_{1},F_{0})(n,x),(\delta ^{\ast }t)(n,x)\in A_{1}$ $\simeq Kers{^{\prime }}\rtimes
A_{0}^{^{\prime }}$ such that $t(F_{1},F_{0})(n,x)=s(\delta ^{\ast }t)(n,x),$ we have $%
(F_{1},F_{0})(n,x)\circ ^{^{\prime }}\delta ^{\ast }t(n,x)$ $=(F_{1}(n)+\delta
t(n),F_{0}(x))$ and $-(F_{1},F_{0})(n,x)=(-F_{1}(n),t^{^{\prime
}}F_{1}(n)+F_{0}(x))$ and then, since
\begin{equation*}
(F_{1},F_{0})(n,x)\circ ^{^{\prime }}\delta ^{\ast }t(n,x)=\delta ^{\ast
}s(n,x)\circ ^{^{\prime }}(G_{1},G_{0})(n,x)
\end{equation*}%
we have
\begin{equation*}
\begin{array}{lll}
\delta ^{\ast }t(n,x) & = & -(F_{1},F_{0})(n,x)\circ ^{^{\prime }}\delta ^{\ast
}s(n,x)\circ ^{^{\prime }}(G_{1},G_{0})(n,x) \\
& = & (-F_{1}(n)+h(x)+G_{1}(n),t^{^{\prime }}F_{1}(n)+F_{0}(x))%
\end{array}%
\end{equation*}%
and
\begin{equation*}
-e^{^{\prime }}F_{0}t(n,x)=(0,t^{^{\prime }}f_{1}(n)+f_{0}(x)).
\end{equation*}%
Hence we get%
\begin{equation*}
\begin{array}{lll}
\delta ^{\ast }t(n,x)-e^{^{\prime }}F_{0}t(n,x) & = & (I_{t^{^{\prime
}}F_{1}(n)+F_{0}(x)}\circ \delta t)(n,x) \\
& = & \delta ^{\ast }t(n,x).%
\end{array}%
\end{equation*}%
Then%
\begin{equation*}
\begin{array}{lll}
h^{\ast }(t(n,x)) & = & \delta ^{\ast }(t(n,x))-e^{^{\prime }}(s^{^{\prime
}}\delta ^{\ast })(t(n,x)) \\
& = & \delta ^{\ast }t(n,x)-e^{^{\prime }}F_{0}t^{\ast }(n,x) \\
& = & \delta ^{\ast }t(n,x) \\
& = & (-F_{1}(n)+h(x)+G_{1}(n),t^{^{\prime }}F_{1}(n)+F_{0}(x))\text{ \ \ \
\ \ \ \ \ \ \ \ \ \ \ \ \ \ \ \ (1)}%
\end{array}%
\end{equation*}%
and%
\begin{equation*}
\begin{array}{lll}
h^{\ast }(t(n,x)) & = & h^{\ast }(\partial (n)+x)) \\
& = & (h(\partial (n)+x)),f_{0}(\partial (n)+x)) \\
& = & (h(\partial (n))+h(x),f_{0}(\partial (n))+f_{0}(x)) \\
& = & (h(\partial (n))+h(x),t^{^{\prime }}F_{1}(n)+F_{0}(x)).\text{ \ \ \ \ \
\ \ \ \ \ \ \ \ \ \ \ \ \ \ (2)}%
\end{array}%
\end{equation*}%
Therefore from (1) and (2) we have%
\begin{equation*}
h(\partial (n))+h(x)=-F_{1}(n)+h(x)+G_{1}(n)
\end{equation*}%
and
\begin{equation*}
h(\partial (n))=-F_{1}(n)+G_{1}(n).
\end{equation*}%
Then
\begin{equation*}
g_{1}(n)=f_{1}(n)+h\partial (n).
\end{equation*}%
Hence%
\begin{equation*}
\begin{array}{cccc}
h: & A_{0} & \longrightarrow  & Kers^{^{\prime }} \\
& x & \longmapsto  & h(x)=\delta (x)-e^{^{\prime }}F_{0}(x)%
\end{array}%
\end{equation*}%
is a homotopy connecting $f=(f_{1},f_{0}):(Kers\overset{\partial }{%
\longrightarrow }A_{0})\longrightarrow (Kers^{^{\prime }}\overset{\partial
^{^{\prime }}}{\longrightarrow }A_{0}^{^{\prime }})$ to $%
g=(g_{1},g_{0}):(Kers\overset{\partial }{\longrightarrow }%
A_{0})\longrightarrow (Kers^{^{\prime }}\overset{\partial ^{^{\prime }}}{%
\longrightarrow }A_{0}^{^{\prime }}).$

Let $F\overset{\delta }{\longrightarrow }G$ and $G\overset{\delta ^{^{\prime }}}{%
\longrightarrow }H.$ Then we have
\begin{equation*}
\begin{array}{lll}
\Gamma (\delta \ast \delta ^{^{\prime }})(x) & = & (\delta \ast \delta ^{^{\prime
}})(x)-e^{^{\prime }}(s^{^{\prime }}\delta \ast \delta ^{^{\prime }} )(x) \\
& = & \delta (x)+\delta ^{^{\prime }}(x)-e^{^{\prime }}(t^{^{\prime }}\delta
)(x)-e^{^{\prime }}(s^{^{\prime }}\delta )(x) \\
& = & \delta (x)+\delta ^{^{\prime }}(x)-e^{^{\prime }}(s^{^{\prime }}\delta
^{^{\prime }})(x)-e^{^{\prime }}(s^{^{\prime }}\delta )(x) \\
& = & (\delta (x)-e^{^{\prime }}(s^{^{\prime }}\delta )(x))+(\delta ^{^{\prime
}}(x)-e^{^{\prime }}(s^{^{\prime }}\delta ^{^{\prime }})(x)) \\
& = & \Gamma (\delta )(x)+\Gamma (\delta ^{^{\prime }})(x)%
\end{array}%
\end{equation*}%
for all $x\in A_{0}.$
\end{pf}

\begin{thm}
Let $\Psi :XMod_{k}\longrightarrow 2A\lg $ be the functor as mentioned in
Theorem 3.1 and $h$ be homotopy connecting $f:(G,C,\partial )\longrightarrow
(G^{^{\prime }},C^{^{\prime }},\partial ^{^{\prime }})$ to $g:(G,C,\partial
)\longrightarrow (G^{^{\prime }},C^{^{\prime }},\partial ^{^{\prime }})$.
Then
\begin{equation*}
\begin{array}{ccccc}
\Psi (h)=\delta & : & C & \longrightarrow & G^{^{\prime }}\rtimes C^{^{\prime }}
\\
&  & x & \longmapsto & \delta (x)=(h(x),f_{0}(x))%
\end{array}%
\end{equation*}%
is a homotopy of corresponding 2-algebra morphisms.
\end{thm}

\begin{pf}
We first show that $\delta $ is an algebra morphism. For $x,x^{^{\prime }}\in C$%
\begin{equation*}
\begin{array}{lll}
\delta (xx^{^{\prime }}) & = & (h(xx^{^{\prime }}),f_{0}(xx^{^{\prime }})) \\
& = & (f_{0}(x)\blacktriangleright h(x^{^{\prime }})+f_{0}(x^{^{\prime
}})\blacktriangleright h(x)+h(x)h(x^{^{\prime }}),~f_{0}(x)f_{0}(x^{^{\prime
}})) \\
& = & (h(x),f_{0}(x))(h(x^{^{\prime }}),f_{0}(x^{^{\prime }})) \\
& = & \delta (x)\delta (x^{^{\prime }}).%
\end{array}%
\end{equation*}%
Now we show that

1) $s^{^{\prime }}\delta =F_{0}\qquad $2) $t^{^{\prime }}$ $\delta
=G_{0}\qquad $3) $(f_{1},f_{0})\circ ^{\prime }\delta t=\delta s\circ ^{\prime }(g_{1},g_{0})$

1)For all $x\in C$,
\begin{equation*}
\begin{array}{lll}
s^{^{\prime }}\delta (x) & = & s^{^{\prime }}(h(x),f_{0}(x)) \\
& = & f_{0}(x)=F_{0}(x),%
\end{array}%
\end{equation*}%
2)For all $x\in C, $
\begin{equation*}
\begin{array}{lll}
t^{^{\prime }}\delta (x) & = & t^{^{\prime }}(h(x),f_{0}(x)) \\
& = & t^{^{\prime }}(h(x))+f_{0}(x) \\
& = & \partial ^{^{\prime }}h(x)+f_{0}(x) \\
& = & g_{0}(x)=G_{0}(x),%
\end{array}%
\end{equation*}%
3)For all $x\in C,a\in G$, since $t^{^{\prime }}(f_{1}(a),f_{0}(x))=\partial ^{^{\prime
}}f_{1}(a)+f_{0}(x)$,
\begin{equation*}
\begin{array}{lll}
s^{^{\prime }}(\delta t(a,x)) & = & s^{^{\prime }}(\delta (\partial (a)+x)) \\
& = & s^{^{\prime }}(h(\partial (a)+x),f_{0}(\partial (a)+x)) \\
& = & f_{0}(\partial (a)+x) \\
& = & f_{0}(\partial (a))+f_{0}(x) \\
& = & \partial ^{^{\prime }}f_{1}(a)+f_{0}(x)%
\end{array}%
\end{equation*}%
then $t^{^{\prime }}(f_{1}(a),f_{0}(x))=s^{^{\prime }}(\delta t(a,x))$ and $%
(f_{1},f_{0})$ , $\delta t$ are composable pairs. Also since
\begin{equation*}
\begin{array}{lll}
t^{^{\prime }}(\delta s(a,x)) & = & t^{^{\prime }}(\delta (x))=t^{^{\prime
}}(h(x),f_{0}(x)) \\
& = & \partial ^{^{\prime }}(h(x))+f_{0}(x) \\
& = & g_{0}(x)%
\end{array}%
\end{equation*}%
and $s^{^{\prime }}(g_{1}(a),g_{0}(x))=g_{0}(x)$ then $t^{^{\prime }}(\delta
s)=s^{^{\prime }}(g_{1},g_{0})$ and $\delta s,(g_{1},g_{0})$ are composable
pairs.

Therefore we get%
\begin{equation*}
\begin{array}{ccc}
(f_{1}(a),f_{0}(x))\circ ^{^{\prime }}\delta t(a,x) & = & (f_{1}(a)+h(\partial
(a)+x),f_{0}(x))%
\end{array}%
\end{equation*}%
and
\begin{equation*}
\begin{array}{ccc}
\delta s(a,x)\circ ^{^{\prime }}(g_{1}(a),g_{0}(x)) & = & (f_{1}(a)+h(\partial
(a)+x),f_{0}(x)).%
\end{array}%
\end{equation*}%
Then $(f_{1},f_{0})\circ ^{^{\prime }}\delta t=\delta s\circ ^{^{\prime
}}(g_{1},g_{0}).$ So%
\begin{equation*}
\begin{array}{cccc}
\delta : & C & \longrightarrow & G^{^{\prime }}\rtimes C^{^{\prime }} \\
& c & \longmapsto & \delta (x)=(h(x),f_{0}(x))%
\end{array}%
\end{equation*}%
is a homotopy connecting $F=((f_{1},f_{0}),f_{0})$ to $%
G=((g_{1},g_{0}),g_{0}).$

Let $f\overset{h}{\longrightarrow }g$ and $g\overset{h^{^{\prime }}}{%
\longrightarrow }u.$ Then we have
\begin{equation*}
\begin{array}{lll}
\Psi (h+h^{^{\prime }})(x) & = & ((h+h^{^{\prime }})(x),f_{0}(x)) \\
& = & (h(x)+h^{^{\prime }}(x),f_{0}(x)) \\
& = & (h(x),f_{0}(x))+(h^{^{\prime }}(x),g_{0}(x))-(0,g_{0}(x)) \\
& = & \Psi (h)(x)+\Psi (h^{^{\prime }})(x)-e^{^{\prime }}(t^{^{\prime
}}(\Psi )(h))(x) \\
& = & (\Psi (h)\ast \Psi (h))(x).%
\end{array}%
\end{equation*}
\end{pf}

\bigskip

\bigskip

\bigskip \.{I}brahim \.{I}lker Ak\c{c}a \newline
\textit{Department of Mathematics-Computer,\newline
Faculty of Science and Letters, \newline
Eski\c{s}ehir Osmangazi University, \newline
26480, Eski\c{s}ehir, Turkey } \newline
e-Mail: iakca@ogu.edu.tr, \newline

\bigskip

Ummahan Ege Arslan \newline
\textit{Department of Mathematics-Computer,\newline
Faculty of Science and Letters, \newline
Eski\c{s}ehir Osmangazi University, \newline
26480, Eski\c{s}ehir, Turkey } \newline
e-Mail: uege@ogu.edu.tr. \newline

\end{document}